\documentclass[11pt]{article}
\usepackage{amsmath,amsfonts}
\usepackage{verbatim}
\usepackage{latexsym}
\usepackage{graphicx}
\usepackage{float}
\usepackage{color}
\usepackage{cite}
\allowdisplaybreaks
\makeatletter
\@addtoreset{equation}{section}
\renewcommand{\thesection}{\arabic{section}}

\begin{document}

\newcommand{\E}{\mathbb{E}}
\newcommand{\PP}{\mathbb{P}}
\newcommand{\RR}{\mathbb{R}}

\newtheorem{theorem}{Theorem}[section]
\newtheorem{lemma}[theorem]{Lemma}
\newtheorem{coro}[theorem]{Corollary}
\newtheorem{defn}[theorem]{Definition}
\newtheorem{assp}[theorem]{Assumption}
\newtheorem{expl}[theorem]{Example}
\newtheorem{prop}[theorem]{Proposition}
\newtheorem{rmk}[theorem]{Remark}

\newcommand\tq{{\scriptstyle{3\over 4 }\scriptstyle}}
\newcommand\qua{{\scriptstyle{1\over 4 }\scriptstyle}}
\newcommand\hf{{\textstyle{1\over 2 }\displaystyle}}
\newcommand\hhf{{\scriptstyle{1\over 2 }\scriptstyle}}

\newcommand{\proof}{\noindent {\it Proof}. }
\newcommand{\eproof}{\hfill $\Box$} 

\def\a{\alpha} \def\g{\gamma}
\def\e{\varepsilon} \def\z{\zeta} \def\y{\eta} \def\o{\theta}
\def\vo{\vartheta} \def\k{\kappa} \def\l{\lambda} \def\m{\mu} \def\n{\nu}
\def\x{\xi}  \def\r{\rho} \def\s{\sigma}
\def\p{\phi} \def\f{\varphi}   \def\w{\omega}
\def\q{\surd} \def\i{\bot} \def\h{\forall} \def\j{\emptyset}

\def\be{\beta} \def\de{\delta} \def\up{\upsilon} \def\eq{\equiv}
\def\ve{\vee} \def\we{\wedge}

\def\F{{\cal F}}
\def\T{\tau} \def\G{\Gamma}  \def\D{\Delta} \def\O{\Theta} \def\L{\Lambda}
\def\X{\Xi} \def\S{\Sigma} \def\W{\Omega}
\def\M{\partial} \def\N{\nabla} \def\Ex{\exists} \def\K{\times}
\def\V{\bigvee} \def\U{\bigwedge}

\def\1{\oslash} \def\2{\oplus} \def\3{\otimes} \def\4{\ominus}
\def\5{\circ} \def\6{\odot} \def\7{\backslash} \def\8{\infty}
\def\9{\bigcap} \def\0{\bigcup} \def\+{\pm} \def\-{\mp}
\def\la{\langle} \def\ra{\rangle}

\def\tl{\tilde}
\def\trace{\hbox{\rm trace}}
\def\diag{\hbox{\rm diag}}
\def\for{\quad\hbox{for }}
\def\refer{\hangindent=0.3in\hangafter=1}

\newcommand\wD{\widehat{\D}}
\title{
\bf The randomized Milstein scheme for stochastic Volterra integral equations with weakly singular kernels
}

\author{
{\bf Zhaohang Wang${}^{1}$,  Zhuoqi Liu${}^{2,}$, Shuaibin Gao${}^2$, Junhao Hu${}^1$\thanks{The corresponding author. Email: junhaohu74@163.com} }
\\
${}^1$ School of Mathematics and Statistics, \\
South-Central Minzu University, \\
Wuhan, 430074, China.\\
${}^2$ Department of Mathematics, \\
Shanghai Normal University, \\
Shanghai, 200234, China. \\
 }

\date{}

\maketitle

\begin{abstract}
This paper focuses on the randomized Milstein scheme for approximating solutions to stochastic Volterra integral equations with weakly singular kernels, where the drift coefficients are non-differentiable. An essential component of the error analysis involves the utilization of randomized quadrature rules for stochastic integrals to avoid the Taylor expansion in drift coefficient functions. Finally, we implement the simulation of multiple singular stochastic integral in the numerical experiment by applying the Riemann-Stieltjes integral.

 \par \noindent

\end{abstract}

\section{Introduction}

Over the past several decades, stochastic Volterra integral equations (SVIEs) have been widely applied in numerous branches of science such as  mathematical finance, physics, biology, engineering and so on \cite{H07, PJ98}. However, the explicit solution of SVIE is difficult to obtain due to the nonlinear property, so it is necessary to develop numerical schemes to approximate the explicit solution \cite{G19, K21, L17, M15, W22}. Especially,  the numerical solutions of stochastic Volterra integral equations with weakly singular kernel (SVIEwWSKs) have been studied by many researchers. The $\theta$-Euler-scheme and the Milstein scheme for SVIEwWSKs were investigated in \cite{L22}. Further,  both the Euler-Maruyama(EM) scheme and the Milstein scheme of a more general form of SVIEs were studied in \cite{R21}.  Later, the fast EM scheme for weakly singular SVIEs with variable exponent was introduced in \cite{L23}. The EM scheme for weakly singular stochastic fractional integro-differential equations was investigated in \cite{D22}. And the fast EM scheme for SVIEs with singular and Hölder continuous kernels was studied in \cite{W23}. For further references about SVIEs, please refer to \cite{Z08, Y19, Z20, ZL20, X18, W11, W17}

In addition, certain randomized Euler and Runge-Kutta scheme have been studied for deterministic differential equation\cite{D11, HM08, JN09, K06, KW17, S90, S95}. The randomized EM scheme for scalar stochastic differential equations (SDEs) with Carathéodory type  drift coefficient functions was investigated in \cite{P14}. Futhermore, the randomized EM schemes for scalar SDEs with drift coefficient which is Lipschitz continuous with respect to the space variable but only measurable with respect to the time variable were introduced in \cite{P15, PP15}. And the randomized Euler scheme of SDEs for which the drift and diffusion coefficients are perturbed by some deterministic noise was studied in \cite{MP17}. 

In this paper, inspired by the drift-randomized Milstein scheme introduced in \cite{KW19}, we apply the randomized scheme to the classical Milstein scheme resulting in the randomized Milstein scheme for SVIEwWSK. Compared with the differentiable condition for the drift coefficient function in \cite{L22}, in this paper the function $b$ is not necessarily differentiable which can also be observed in the numerical simulation. Note that the drift and diffusion coefficients in \cite{KW19} are temporal Hölder continuous, hence it is more challenging to cope with the singular kernels in the drift and diffusion coefficient $(t-s)^{-\a}$ and $(t-s)^{-\be}$ in our paper as they will tend to infinity as $s$ tends to $t$. The main result of this paper shows that, under Assumptions 2.3 and 2.4, the convergence rate of the randomized Milstein scheme for SVIEsWSK can not exceed $\min\{1-2\be, 1-\a\}$. Moreover, testing the convergence rate still remains a problem in \cite{L22, LL23} since the simulation of the singular stochastic integral $\int_{t_{j-1}}^{t_j}(t_n-s)^{-\be}\s^{'}(X_{j-1}^{h})\Big(\sum_{k=1}^{j-1}\int_{t_{k-1}}^{t_k}[(s-r)^{-\be}-(t_{j-1}-r)^{-\be}]\s(X_{k-1}^{h})dBr\Big)dB_s$ is a significant challenge. However, this problem is addressed in our paper by using the Riemann-Stieltjes integral which we will discuss in more details in Section 4 through a numerical experiment. 

The rest of the paper is structured as follows. In Section 2, we introduce some notations, assumptions and a few important lemmas that are useful in the proof later. Section 3 aims to get the final convergence rate of the randomized Milstein scheme. In section 4, a numerical example will be conducted to validate the effectiveness of the theoretical results.

\section{Preliminaries}

If $A$ is a vector or matrix, its transpose is denoted by $A^T$.
If $x\in R^n$, then $|x|$ is its Euclidean norm.
If $A$ is a matrix, let $|A| = \sqrt{\trace(A^TA)}$ be its trace norm.
If $a,b$ are both real numbers, then $a\ve b:=\max\{a,b\}$, $a\we b:=\min\{a,b\}$.
Let $\mathbb{S}_n=\{1,\cdots, n\}$, $\mathbb{S}^0_n=\{0, 1,\cdots, n\}$ for any $n\in\mathbb{N}$. If $s\in\RR$, then $\lfloor s\rfloor$ is the greatest integer no more than $s$. Denote by $C^2(\RR^d)$ the family of twice continuously differentiable functions in $\RR^d$.
For $p>0$ and $t\ge 0$, denote by $L^p(\W, \cal F, \PP)=$ $L^p(\W)$ the family of
 $R^n$-valued random variables $X$ such that
$$
\|X\|_{L^p(\W)}:=(\E|X|^p)^{\frac{1}{p}}=\Big(\int_\W|X(w)|^pd\PP(w)\Big)^\frac{1}{p}<\8.
$$

In this paper, consider the following SVIEwWSK
\begin{equation} \label{1.1}
X(t) =X_0+\int_0^t(t-s)^{-\a} b(X(s))ds + \int_0^t(t-s)^{-\beta} \s(X(s))dB_s, \quad t\in[0,T],
\end{equation}
where $\{B_t\}_{t\in[0,T]}$ is an m-dimensional standard Brownian motion defined on the complete filtered probability space
$(\W_B ,{\cal F}^B, \{{\cal F}_t^B\}_{t\in[0,T]}, \PP_B)$. And $b:\RR^d\to\RR^d$ and $\s:\RR^d\to\RR^{d\times m}$ are Borel measurable functions. Let $\a, \be\in(0,\frac{1}{2})$ be two given positive real numbers and $X_0$ denote the initial value satisfying $\E_B|X_0|^2<\8$.

To begin the introduction of the randomized Milstein scheme in this paper, firstly we partition the interval $[0,T]$ into $N$ equidistant intervals with stepsize $h=\frac{T}{N}$, i.e.,
\begin{equation}\label{I}
I_h:=\{t_n:=nh, n=0,1,\cdots,N\}.
\end{equation}
Further, let$\{\tau_j\}_{j\in\mathbb{N}}$ be an  family of independent identically distributed ${\cal U}(0,1)$-distributed random variables defined on another filtered  probability space $(\W_\T ,{\cal F}^\T, \{{\cal F}_j^\T\}_{j\in \mathbb{N}}, \PP_\T)$, where ${\cal U}(0,1)$ presents the uniform distribution on the interval $(0,1)$ and ${\cal F}_j^\T$ is the $\s$-algebra generated by $\{\T_1,\cdots,\T_j\}$. The random variables $\{\tau_j\}_{j\in\mathbb{N}}$ stand for the artificially added random input for the new method, which we suppose to be independent of the randomness present in (\ref{1.1}).

Then the stochastic process yielded from the numerical scheme will be defined on the product probability space
\begin{equation}\label{1.2}
(\W ,{\cal F},\PP):=(\W_B \K\W_\T,{\cal F}^B\otimes{\cal F}^\T,\PP_B\otimes\PP_\T)
\end{equation}
Moreover, for each partition $I_h$, denote the associated discrete-time filtration$ \{{\cal F}_n^h\}_{n\in \mathbb{S}_N }$ on $(\W ,{\cal F},\PP)$ by
$$
{\cal F}_n^h:={\cal F}_{t_n}^B\otimes{\cal F}_n^\T,\quad  n\in \mathbb{S}^0_N
$$
Then, the randomized Milstein scheme on the partition is given by
\begin{align} \label{1.3}
X_n^h&=X_0+\sum_{j=1}^nh(t_n-(t_{j-1}+\T_jh))^{-\a}b(Y_j^{h,\T})+\sum_{j=1}^n\int_{t_{j-1}}^{t_j}(t_n-s)^{-\be}\s(X_{j-1}^{h})dB_s\nonumber\\
&+\sum_{j=1}^n\int_{t_{j-1}}^{t_j}(t_n-s)^{-\be}\s^{'}(X_{j-1}^{h})\Big(\sum_{k=1}^{j-1}\int_{t_{k-1}}^{t_k}[(s-r)^{-\a}-(t_{j-1}-r)^{-\a}]b(X_{k-1}^{h})dr\Big)dB_s\nonumber\\
&+\sum_{j=1}^n\int_{t_{j-1}}^{t_j}(t_n-s)^{-\be}\s^{'}(X_{j-1}^{h})\Big(\int_{t_{j-1}}^s(s-r)^{-\a}b(X_{j-1}^{h})dr\Big)dB_s\nonumber\\
&+\sum_{j=1}^n\int_{t_{j-1}}^{t_j}(t_n-s)^{-\be}\s^{'}(X_{j-1}^{h})\Big(\sum_{k=1}^{j-1}\int_{t_{k-1}}^{t_k}[(s-r)^{-\be}-(t_{j-1}-r)^{-\be}]\s(X_{k-1}^{h})dBr\Big)dB_s\nonumber\\
&+\sum_{j=1}^n\int_{t_{j-1}}^{t_j}(t_n-s)^{-\be}\s^{'}(X_{j-1}^{h})\Big(\int_{t_{j-1}}^s(s-r)^{-\be}\s(X_{j-1}^{h})dBr\Big)dB_s,\nonumber\\
\text{and}\nonumber\\
Y_j^{h,\T}&=X_{j-1}^h+\sum_{i=1}^{j-1}\int_{t_{i-1}}^{t_i}[(t_{j-1}+\T_jh-s)^{-\a}-(t_{j-1}-s)^{-\a}]b(X_{i-1}^h)ds\nonumber\\
&+\sum_{i=1}^{j-1}\int_{t_{i-1}}^{t_i}[(t_{j-1}+\T_jh-s)^{-\be}-(t_{j-1}-s)^{-\be}]\s(X_{i-1}^h)dB_s\nonumber\\
&+\int_{t_{j-1}}^{t_{j-1}+\T_jh}(t_{j-1}+\T_jh-s)^{-\a}b(X_{j-1}^h)ds+\int_{t_{j-1}}^{t_{j-1}+\T_jh}(t_{j-1}+\T_jh-s)^{-\be}\s(X_{j-1}^h)dB_s,
\end{align}
where the initial value $X_0^h=X_0$ and $n\in\mathbb{S}_N$.

Let $L^p(\W_B, {\cal F}^B and \PP_B)=L^p(\W_B)$, $L^p(\W_\T, {\cal F}^\T, \PP_\T)=L^p(\W_\T)$. In the case of the product probability space$(\W ,{\cal F},\PP)$ introduced in (\ref{1.2}), an application of Funini's theorem shows that
$$
\E[X]=\E_B[\E_\T[X]]=\E_\T[\E_B[X]],\quad X\in L^p(\W),
$$
where $\E_\T$, $\E_B$ are the expectation with respect to $\PP_\T$, $\PP_B$, respectively.

The following discrete Burkholder-Davis-Gundy (BDG) inequality is an important tool  \cite{B66}.
\begin{theorem}
 For every discrete time martingale $(X_n)_{n\in\mathbb{N}}$ and for every $n\in\mathbb{N}$ we have
$$
c_p\|[X]_n^{\frac{1}{2}}\|_{L^p(\W)}\le\|\max_{j\in\{0,\cdots,n\}}|X^j|\|_{L^p(\W)}\le C_p\|[X]_n^{\frac{1}{2}}\|_{L^p(\W)},
$$
where $p\in(1,\8)$, $c_p$ and $C_p$ are positive constants as well as $[X]_n=|X^0|^2+\sum_{k=1}^n|X^k-X^{k-1}|^2$ is the quadratic variation of $(X_n)_{n\in\mathbb{N}}$.
\end{theorem}

 The following generalized discrete type of Gronwall inequality derived from \cite{B04} is needed to cope with the singular kernels.
\begin{lemma}
Assume that $a>0, 0<\g<1$ are two positive numbers, and $\{B_n\}, \{b_n\}$ are two non-negative sequences satisfying the following inequality
$$
B_n\le b_n+ah^{1-\g}\sum_{j=1}^{n}(n+1-j)^{-\g}B_{j-1},
$$
for all $0\le n\le N$. Then
$$
B_n\le E_{1-\g}(\G(1-\g)(nh)^{1-\g}a)b_n,
$$
where $\G(a)=\int_0^\8e^{-s}s^{m-1}ds, a>0$ is the Euler-Gamma function and
$$
E_a(x)=:\sum_{k=0}^\8\frac{1}{\G(mk+1)}x^k,\quad m>0,
$$
is the Mittag-Leffler function of $x\in\RR$
\end{lemma}

\begin{assp}\label{a1}
 There exists a positive constant $\hat{L}_1$ such that, for all $x,y\in\RR^d$,
\begin{equation}
|b(x)-b(y)|\ve|\s(x)-\s(y)|\le\hat{L}_1|x-y|.\nonumber
\end{equation}
\end{assp}

\begin{assp}
The diffusion coefficient $\s$ is in $C^2(\RR^d)$. Moreover, there exists a positive constant $\hat{L}_2$ such that, for all $x\in\RR^d$, \begin{align*}
&|\nabla \s(x)|\le \hat{L}_2,\\
&|\nabla \s(x)-\nabla \s(y)|\le \hat{L}_2|x-y|,
\end{align*}

\end{assp}
Using (A1), one has $|b(x)|\ve|\s(x)|\le C_0(1+|x|)$ for some constant $C_0>0$.

The existence, uniqueness, boundness and Hölder continuity of the solution in \cite{L22} will  be stated in the following theorem. 
\begin{theorem}
Let Assumption 2.3 hold, then there exists a unique solution $X(t)$ to (\ref{1.1}), and the solution has the following properties
\begin{align}
&\sup_{t\in[0,T]}\E_B|X(t)|^2\le C_1,\\
&\|X(t)-X(s)\|_{L^2(\W_B)}\le C_2|t-s|^{\g},\quad \forall  0\le s< t\le T ,
\end{align}
where $\g=\min\{\frac{1}{2}-\be,1-\a\}$ and $C_1, C_2$ are two positive real constant numbers.
\end{theorem}

Let $I_h$ be the same partition as in (\ref{I}), then we define the increment function
$\Phi_n^{h,\T}$ of the n-th step for each $n\in\mathbb{S}_N$ as following
\begin{align}\label{2.1}
\Phi_n^{h,\T}(y)&=\sum_{j=1}^nh(t_n-(t_{j-1}+\T_jh))^{-\a}b(\Psi_j^{h,\T}(y))+\sum_{j=1}^n\int_{t_{j-1}}^{t_j}(t_n-s)^{-\be}\s(y_{j-1})dB_s\nonumber\\
&+\sum_{j=1}^n\int_{t_{j-1}}^{t_j}(t_n-s)^{-\be}\s^{'}(y_{j-1})\Big(\sum_{k=1}^{j-1}\int_{t_{k-1}}^{t_k}[(s-r)^{-\a}-(t_{j-1}-r)^{-\a}]b(y_{k-1})dr\Big)dB_s\nonumber\\
&+\sum_{j=1}^n\int_{t_{j-1}}^{t_j}(t_n-s)^{-\be}\s^{'}(y_{j-1})\Big(\int_{t_{j-1}}^s(s-r)^{-\a}b(y_{j-1})dr\Big)dB_s\nonumber\\
&+\sum_{j=1}^n\int_{t_{j-1}}^{t_j}(t_n-s)^{-\be}\s^{'}(y_{j-1})\Big(\sum_{k=1}^{j-1}\int_{t_{k-1}}^{t_k}[(s-r)^{-\be}-(t_{j-1}-r)^{-\be}]\s(y_{k-1})dBr\Big)dB_s\nonumber\\
&+\sum_{j=1}^n\int_{t_{j-1}}^{t_j}(t_n-s)^{-\be}\s^{'}(y_{j-1})\Big(\int_{t_{j-1}}^s(s-r)^{-\be}\s(y_{j-1})dBr\Big)dB_s,
\end{align}
and
\begin{align}
\Psi_j^{h,\T}(y)&=y_{j-1}+\sum_{i=1}^{j-1}\int_{t_{i-1}}^{t_i}[(t_{j-1}+\T_jh-s)^{-\a}-(t_{j-1}-s)^{-\a}]b(y_{i-1})ds\nonumber\\
&+\sum_{i=1}^{j-1}\int_{t_{i-1}}^{t_i}[(t_{j-1}+\T_jh-s)^{-\be}-(t_{j-1}-s)^{-\be}]\s(y_{i-1})dB_s\nonumber\\
&+\int_{t_{j-1}}^{t_{j-1}+\T_jh}(t_{j-1}+\T_jh-s)^{-\a}b(y_{j-1})ds+\int_{t_{j-1}}^{t_{j-1}+\T_jh}(t_{j-1}+\T_jh-s)^{-\be}\s(y_{j-1})dB_s,
\end{align}
where $y=(y_0,\cdots,y_{n-1})\in\RR^{d\times n}$, $\T_j\in(0,1)$ for $j\in\mathbb{S}_n$. Then (\ref{1.3}) can be rewritten by
\begin{equation}
\left\{
\begin{array}{ll}
X_n^h=X_0+\Phi_n^{h,\T}(z),\quad z=(X_0^h,\cdots,X_{n-1}^h), n\in\mathbb{S}_N,\\
X_0^h=X_0.
\end{array}
\right.
\end{equation}

The following lemma ensures that $\{X_n^h\}_{n\in\mathbb{S}_N}$ is an adapted sequence in $L^p(\W)$.
\begin{lemma}
Let Assumption 2.3 and Assumption 2.4 hold. If $y=(y_0,\cdots,y_{n-1})$ satisfies $y_j\in L^2(\W ,{\cal F}_j^h,\PP)$ for all $j\in\mathbb{S}^0_{n-1}$, then
$$
\Phi_n^{h,\T}(y)\in L^2(\W ,{\cal F}_n^h,\PP)
$$
\end{lemma}
\noindent
{\it Proof}. We prove $\Psi_j^{h,\T}(y)\in L^2(\W ,{\cal F}_{j}^h,\PP)$ at first. It is easy to see that $\Psi_j^{h,\T}(y)$ is ${\cal F}_j^h$-measurable, thus it remains to prove the $L^2$-boundedness of $\Psi_j^{h,\T}(y)$. Note that
\begin{align}
\|\Psi_j^{h,\T}(y)\|_{L^2(\W)}&\le\|y_{j-1}\|_{L^2(\W)}
+\bigg{\|}\int_0^{t_{j-1}}[(t_{j-1}+\T_jh-s)^{-\a}-(t_{j-1}-s)^{-\a}]b(y_{\lfloor s/h\rfloor})ds\bigg{\|}_{L^2(\W)}\nonumber\\
&+\bigg{\|}\int_0^{t_{j-1}}[(t_{j-1}+\T_jh-s)^{-\a}-(t_{j-1}-s)^{-\a}]\s(y_{\lfloor s/h\rfloor})dB_s\bigg{\|}_{L^2(\W)}\nonumber\\
&+\bigg{\|}\int_{t_{j-1}}^{t_{j-1}+\T_jh}(t_{j-1}+\T_jh-s)^{-\a}b(y_{j-1})ds\bigg{\|}_{L^2(\W)}\nonumber\\
&+\bigg{\|}\int_{t_{j-1}}^{t_{j-1}+\T_jh}(t_{j-1}+\T_jh-s)^{-\be}\s(y_{j-1})dB_s\bigg{\|}_{L^2(\W)}\nonumber\\
&=:\|y_{j-1}\|_{L^2(\W)}+I_1^1+I_2^1+I_3^1+I_4^1\nonumber.
\end{align}
First, for the estimate of $I_1^1$ we have
\begin{align}
I_1^1\le C_0\bigg(\E\Big(\int_0^{t_{j-1}}|(t_{j}-s)^{-\a}-(t_{j-1}-s)^{-\a}|(1+|y_{\lfloor s/h\rfloor}|)ds\Big)^2\bigg)^{\frac{1}{2}}.
\end{align}
Denote
\begin{equation}
\rho_{t_j,t_{j-1}}=\int_0^{t_{j-1}}|(t_{j}-s)^{-\a}-(t_{j-1}-s)^{-\a}|ds\nonumber.
\end{equation}
Since $\phi(x)=x^2, x>0$ is a convex function, applying Jensen's inequality yields
 \begin{align}
 &\bigg(\frac{1}{\rho_{t_j,t_{j-1}}}\int_0^{t_{j-1}}|(t_{j}-s)^{-\a}-(t_{j-1}-s)^{-\a}|(1+|y_{\lfloor s/h\rfloor}|)ds\bigg)^2\nonumber\\
 \le&\frac{1}{\rho_{t_j,t_{j-1}}}\int_0^{t_{j-1}}|(t_{j}-s)^{-\a}-(t_{j-1}-s)^{-\a}|(1+|y_{\lfloor s/h\rfloor}|)^2ds\nonumber.
 \end{align}
 We therefore have
 \begin{align}
 I_1^1&\le C_0\bigg(\rho_{t_j,t_{j-1}}\int_0^{t_{j-1}}|(t_{j}-s)^{-\a}-(t_{j-1}-s)^{-\a}|\E(1+|y_{\lfloor s/h\rfloor}|)^2ds\bigg)^{\frac{1}{2}}\nonumber\\
 &\le C\rho_{t_j,t_{j-1}}<\8.
 \end{align}
 where $C$ is a constant. Here the estimate of $\rho_{t_j,t_{j-1}}^i$ comes from
 \begin{align}\label{L5.1}
\rho_{t_j,t_{j-1}}&\le2\int_0^{t_{j-1}}(t_{j-1}-s)^{-\a}ds\nonumber\\
&\le\frac{2t_{j-1}^{1-\a}}{1-\a }\le\frac{2T^{1-\a}}{1-\a },
\end{align}
 In the similar way, we can prove that $I_3^1<\8$. For $I_2^1$,  using the BDG inequality and the Jensen inequality, we have
 \begin{align}\label{L4.1}
 I_2^1&=\bigg(\E_{\T}\Big[\E_B\Big|\int_0^{t_{j-1}}[(t_{j-1}+\T_jh-s)^{-\be}-(t_{j-1}-s)^{-\be}]\s(y_{\lfloor s/h\rfloor})dB_s\Big|^2\Big]\bigg)^{\frac{1}{2}}\nonumber\\
 &\le 2C_0\bigg(\E_{\T}\Big[\E_B\Big(\int_0^{t_{j-1}}|(t_j-s)^{-\be}-(t_{j-1}-s)^{-\be}|^2(1+|y_{\lfloor s/h\rfloor}|)^2ds\Big)\Big]\bigg)^{\frac{1}{2}}\nonumber\\
 &\le 2C_0\bigg(\int_0^{t_{j-1}}|(t_j-s)^{-\be}-(t_{j-1}-s)^{-\be}|^2\E(1+|y_{\lfloor s/h\rfloor}|)^2ds\bigg)^{\frac{1}{2}}\nonumber\\
  &\le C\bigg(\int_0^{t_{j-1}}|(t_j-s)^{-\be}-(t_{j-1}-s)^{-\be}|^2ds\bigg)^\frac{1}{2}\nonumber\\
  &\le C\bigg(\int_0^{t_{j-1}}(t_{j-1}-s)^{-2\be}ds\bigg)^\frac{1}{2} <\8.
 \end{align}
  similarly, we can prove that $I_4^1<\8$. Thus,  $\Psi_j^{h,\T}(y)\in L^2(\W ,{\cal F}_j^h,\PP)$. It remains to prove $\Phi_n^{h,\T}(y)\in L^2(\W ,{\cal F}_n^h,\PP)$. Obviously, $\Phi_n^{h,\T}(y)$ is ${\cal F}_n^h$-measurable, now we need to obtain the $L^2$-boundedness of $\Phi_n^{h,\T}(y)$. Note that
 \begin{align}
&\|\Phi_n^{h,\T}(y)\|^2_{L^2(\W)}\nonumber\\
\le&6\Big\|\sum_{j=1}^nh(t_n-(t_{j-1}+\T_jh))^{-\a}b(\Psi_j^{h,\T}(y))\Big\|^2_{L^2(\W)}\nonumber\\
+&6\bigg{\|}\sum_{j=1}^n\int_{t_{j-1}}^{t_j}(t_n-s)^{-\be}\s(y_{j-1})dB_s\bigg{\|}^2_{L^2(\W)}\nonumber\\
+&6\bigg{\|}\sum_{j=1}^n\int_{t_{j-1}}^{t_j}(t_n-s)^{-\be}\s^{'}(y_{j-1})\Big(\sum_{k=1}^{j-1}\int_{t_{k-1}}^{t_k}[(s-r)^{-\a}-(t_{j-1}-r)^{-\a}]b(y_{k-1})dr\Big)dB_s\bigg{\|}^2_{L^2(\W)}\nonumber\\
+&6\bigg{\|}\sum_{j=1}^n\int_{t_{j-1}}^{t_j}(t_n-s)^{-\be}\s^{'}(y_{j-1})\Big(\int_{t_{j-1}}^s(s-r)^{-\a}b(y_{j-1})dr\Big)dB_s\bigg{\|}^2_{L^2(\W)}\nonumber\\
+&6\bigg{\|}\sum_{j=1}^n\int_{t_{j-1}}^{t_j}(t_n-s)^{-\be}\s^{'}(y_{j-1})\Big(\sum_{k=1}^{j-1}\int_{t_{k-1}}^{t_k}[(s-r)^{-\be}-(t_{j-1}-r)^{-\be}]\s(y_{k-1})dB_r\Big)dB_s\bigg{\|}^2_{L^2(\W)}\nonumber\\
+&6\bigg{\|}\sum_{j=1}^n\int_{t_{j-1}}^{t_j}(t_n-s)^{-\be}\s^{'}(y_{j-1})\Big(\int_{t_{j-1}}^s(s-r)^{-\be}\s(y_{j-1})dB_r\Big)dB_s\bigg{\|}^2_{L^2(\W)}\nonumber\\
=:&6(I_1^2+I_2^2+I_3^2+I_4^2+I_5^2+I_6^2).
\end{align}
For $I_1^2$, given the $L^2$-boundedness of $\Psi_j^{h,\T}(y)$, we have
\begin{align}
I_1^2&=\E\Big|\sum_{j=1}^nh(t_n-(t_{j-1}+\T_jh))^{-\a}b(\Psi_j^{h,\T}(y))\Big|^2\nonumber\\
&\le n\sum_{j=1}^n\E_{\T}\big[\E_B|h(t_n-(t_{j-1}+\T_jh))^{-\a}b(\Psi_j^{h,\T}(y))|^2\big]\nonumber\\
&\le 2C_0n\sum_{j=1}^n\E_{\T}\big[|h(t_n-(t_{j-1}+\T_jh))^{-\a}|^2(1+\E_B|\Psi_j^{h,\T}(y)|^2)\big]\nonumber\\
&\le Cn\sum_{j=1}^n\E_{\T}\big|h(t_n-(t_{j-1}+\T_jh))^{-\a}|^2<\8.\nonumber
\end{align}
Here, due to  $\T_j\sim{\cal U}(0,1)$, the estimate of $\E_{\T}\big|h(t_n-(t_{j-1}+\T_jh))^{-\a}|^2$ comes from
\begin{align}\label{T5.1}
&\E_{\T}\big|h(t_n-(t_{j-1}+\T_jh))^{-\a}|^2\nonumber\\
=&h^2\int_0^1(t_n-(t_{j-1}+hs))^{-2\a}ds\nonumber\\
=&h\int_{t_{j-1}}^{t_j}(t_n-s)^{-2\a}ds\nonumber\\
=&\frac{h^{2(1-\a)}}{1-2\a}[(n+1-j)^{1-2\a}-(n-j)^{1-2\a}].
\end{align}
Then, using Assumption 2.3, the BDG inequality and the Jensen inequality we obtain
\begin{align}\label{L4.4}
I_2^2&=\bigg\|\int_0^{t_n}(t_n-s)^{-\be}\s(y_{\lfloor s/h\rfloor})dB_s\bigg\|^2_{L^2(\W)}\nonumber\\
&\le 4C_0\E\int_0^{t_n}(t_n-s)^{-2\be}(1+|y_{\lfloor s/h\rfloor}|)^2ds\nonumber\\
&\le C\int_0^{t_n}(t_n-s)^{-2\be}ds<\8.
\end{align}
Denote $\eta(s)=t_i$ for $t_i\le s\le t_{i+1}$, therefore $|\eta(s)-s|\le h$ holds for all $s\in[0,T]$. We can rewrite $I_3^2$ in a continuous form
$$
I_3^2=\bigg{\|}\int_0^{t_n}(t_n-s)^{-\be}\s^{'}(y_{\lfloor s/h\rfloor})\Big(\int_0^{\eta(s)}[(s-r)^{-\a}-(\eta(s)-r)^{-\a}]b(y_{\lfloor r/h\rfloor})dr\Big)dB_s\bigg{\|}^2_{L^2(\W)}
$$
Applying Assumption 2.3, the BDG inequality and the Jensen inequality, we have
\begin{align}\label{L4.3}
&I_3^2\nonumber\\
\le&4\hat{L}_2\E\int_0^{t_n}(t_n-s)^{-2\be}\Big(\int_0^{\eta(s)}|(s-r)^{-\a}-(\eta(s)-r)^{-\a}||b(y_{\lfloor r/h\rfloor})|dr\Big)^2ds\nonumber\\
\le&C\int_0^{t_n}(t_n-s)^{-2\be}\Big(\int_0^{\eta(s)}|(s-r)^{-\a}-(\eta(s)-r)^{-\a}|dr\Big)^2ds<\8,
\end{align}
where
 \begin{align}\label{L4.2}
&\int_0^{\eta(s)}|(s-r)^{-\a}-(\eta(s)-r)^{-\a}|dr\nonumber\\
 =&\a\int_0^{\eta(s)}\int_{\eta(s)}^{s}(u-r)^{-\a-1}dudr\nonumber\\
 =&\a\int_{\eta(s)}^{s}\int_0^{\eta(s)}(u-r)^{-\a-1}drdu\nonumber\\
 \le&\int_{\eta(s)}^{s}(u-\eta(s))^{-a}du=\frac{1}{1-\a}h^{1-\a}.
 \end{align}
In the similar way, we can prove that $I_4^2<\8, I_5^2<\8, I_6^2<\8$.  Therefore $\Phi_n^{h,\T}(y)\in L^2(\W ,{\cal F}_n^h,\PP)$. The proof is complete.
\eproof

The following lemma can be directly derived from Lemma 2.6.
\begin{lemma}
Let Assumption 2.3 and Assumption 2.4 hold. Then there exists a constant $C_3$ such that 
$$
\max_{n\in\mathbb{S}^0_N}\E|X_n^h|^2\le C_3.
$$
\end{lemma}

\section{Main results}

In this section, we will obtain the main convergence results of the randomized Milstein scheme for SVIEwWSKs. 
\begin{theorem}
Let Assumption 2.3 and Assumption 2.4 hold. Assume that for each $n\in\mathbb{S}_N$, $y=(y_0,\cdots,y_{n-1})$, $\overline{y}=(\overline{y}_0,\cdots,\overline{y}_{n-1})\in\RR^{d\times n}$ satisfy $y_j,\overline{y}_j\in L^2(\W ,{\cal F}_j^h,\PP)$ for $j\in\mathbb{S}^0_{n-1}$, then it holds true that
\begin{equation}
\|\Phi_n^{h,\T}(y)-\Phi_n^{h,\T}(\overline{y})\|_{L^2(\W)}\le CTh^{\g}\sum_{j=1}^n(n+1-j)^{-(\a\ve2\be)}\|y_{j-1}-\overline{y}_{j-1}\|_{L^2(\W)}
\end{equation}
where $C$ is a positive constant depending on $ \hat{L}_1, \a,\be, \|y_j\|_{L^2(\W)}, \|\overline{y}_j\|_{L^2(\W)}, j\in\mathbb{S}^0_{n-1}$.
\end{theorem}
\noindent
{\it Proof}. It follows from (\ref{2.1}) that
\begin{align}\label{T3.1}
&\|\Phi_n^{h,\T}(y)-\Phi_n^{h,\T}(\overline{y})\|_{L^2(\W)}\nonumber\\
\le&\Big\|\sum_{j=1}^nh(t_n-(t_{j-1}+\T_jh))^{-\a}(b(\Psi_j^{h,\T}(y))-b(\Psi_j^{h,\T}(\overline{y})))\Big\|_{L^2(\W)}\nonumber\\
+&\bigg{\|}\sum_{j=1}^n\int_{t_{j-1}}^{t_j}(t_n-s)^{-\be}(\s(y_{j-1})-\s(\overline{y}_{j-1}))dB_s\bigg{\|}_{L^2(\W)}\nonumber\\
+&\bigg{\|}\sum_{j=1}^n\int_{t_{j-1}}^{t_j}(t_n-s)^{-\be}\s^{'}(y_{j-1})\Big(\sum_{k=1}^{j-1}\int_{t_{k-1}}^{t_k}[(s-r)^{-\a}-(t_{j-1}-r)^{-\a}](b(y_{k-1})-b(\overline{y}_{k-1}))dr\Big)dB_s\bigg{\|}_{L^2(\W)}\nonumber\\
+&\bigg{\|}\sum_{j=1}^n\int_{t_{j-1}}^{t_j}(t_n-s)^{-\be}(\s^{'}(y_{j-1})-\s^{'}(\overline{y}_{j-1}))\Big(\sum_{k=1}^{j-1}\int_{t_{k-1}}^{t_k}[(s-r)^{-\a}-(t_{j-1}-r)^{-\a}]b(\overline{y}_{k-1})dr\Big)dB_s\bigg{\|}_{L^2(\W)}\nonumber\\
+&\bigg{\|}\sum_{j=1}^n\int_{t_{j-1}}^{t_j}(t_n-s)^{-\be}\s^{'}(y_{j-1})\Big(\int_{t_{j-1}}^s(s-r)^{-\a}(b(y_{j-1})-b(\overline{y}_{j-1}))dr\Big)dB_s\bigg{\|}_{L^2(\W)}\nonumber\\
+&\bigg{\|}\sum_{j=1}^n\int_{t_{j-1}}^{t_j}(t_n-s)^{-\be}(\s^{'}(y_{j-1})-\s^{'}(\overline{y}_{j-1}))\Big(\int_{t_{j-1}}^s(s-r)^{-\a}b(\overline{y}_{j-1})dr\Big)dB_s\bigg{\|}_{L^2(\W)}\nonumber\\
+&\bigg{\|}\sum_{j=1}^n\int_{t_{j-1}}^{t_j}(t_n-s)^{-\be}\s^{'}(y_{j-1})\Big(\sum_{k=1}^{j-1}\int_{t_{k-1}}^{t_k}[(s-r)^{-\be}-(t_{j-1}-r)^{-\be}](\s(y_{k-1})-\s(\overline{y}_{k-1}))dB_r\Big)dB_s\bigg{\|}_{L^2(\W)}\nonumber\\
+&\bigg{\|}\sum_{j=1}^n\int_{t_{j-1}}^{t_j}(t_n-s)^{-\be}(\s^{'}(y_{j-1})-\s^{'}(\overline{y}_{j-1}))\Big(\sum_{k=1}^{j-1}\int_{t_{k-1}}^{t_k}[(s-r)^{-\be}-(t_{j-1}-r)^{-\be}]\s(\overline{y}_{k-1})dB_r\Big)dB_s\bigg{\|}_{L^2(\W)}\nonumber\\
+&\bigg{\|}\sum_{j=1}^n\int_{t_{j-1}}^{t_j}(t_n-s)^{-\be}\s^{'}(y_{j-1})\Big(\int_{t_{j-1}}^s(s-r)^{-\be}(\s(y_{j-1})-\s(\overline{y}_{j-1}))dB_r\Big)dB_s\bigg{\|}_{L^2(\W)}\nonumber\\
+&\bigg{\|}\sum_{j=1}^n\int_{t_{j-1}}^{t_j}(t_n-s)^{-\be}(\s^{'}(y_{j-1})-\s^{'}(\overline{y}_{j-1}))\Big(\int_{t_{j-1}}^s(s-r)^{-\be}\s(\overline{y}_{j-1})dB_r\Big)dB_s\bigg{\|}_{L^2(\W)}.
\end{align}
Applying Assumption 2.3 and Jensen' inequality, the estimate of $\|\Psi_j^{h,\T}(y)-\Psi_j^{h,\T}(\overline{y})\|_{L^2(\W_B)}$ is obtained by
\begin{align}
&\|\Psi_j^{h,\T}(y)-\Psi_j^{h,\T}(\overline{y})\|_{L^2(\W_B)}
\le\|y_{j-1}-\overline{y}_{j-1}\|_{L^2(\W_B)}\nonumber\\
+&\sum_{i=1}^{j-1}\bigg{\|}\int_{t_{i-1}}^{t_i}[(t_{j-1}+\T_jh-s)^{-\a}-(t_{j-1}-s)^{-\a}](b(y_{i-1})-b(\overline{y}_{i-1}))ds\bigg{\|}_{L^2(\W_B)}\nonumber\\
+&\sum_{i=1}^{j-1}\bigg{\|}\int_{t_{i-1}}^{t_i}[(t_{j-1}+\T_jh-s)^{-\be}-(t_{j-1}-s)^{-\be}](\s(y_{i-1})-\s(\overline{y}_{i-1}))dB_s\bigg{\|}_{L^2(\W_B)}\nonumber\\
+&\bigg{\|}\int_{t_{j-1}}^{t_{j-1}+\T_jh}(t_{j-1}+\T_jh-s)^{-\a}(b(y_{j-1})-b(\overline{y}_{j-1}))ds\bigg{\|}_{L^2(\W_B)}\nonumber\\
+&\bigg{\|}\int_{t_{j-1}}^{t_{j-1}+\T_jh}(t_{j-1}+\T_jh-s)^{-\be}(\s(y_{j-1})-\s(\overline{y}_{j-1}))dB_s\bigg{\|}_{L^2(\W_B)}\nonumber\\
\le&\|y_{j-1}-\overline{y}_{j-1}\|_{L^2(\W_B)}
+2\hat{L}_1\sum_{i=1}^{j-1}\|y_{i-1}-\overline{y}_{i-1}\|_{L^2(\W_B)}\int_{t_{i-1}}^{t_i}(t_{j-1}-s)^{-\a}ds\nonumber\\
+&4\hat{L}_1\sum_{i=1}^{j-1}\|y_{i-1}-\overline{y}_{i-1}\|_{L^2(\W_B)} \bigg(\int_{t_{i-1}}^{t_i}(t_{j-1}-s)^{-2\be}ds\bigg)^{\frac{1}{2}}\nonumber\\
+&\hat{L}_1\|y_{j-1}-\overline{y}_{j-1}\|_{L^2(\W_B)} \int_{t_{j-1}}^{t_j}(t_j-s)^{-\a}ds\nonumber\\
+&2\hat{L}_1\|y_{j-1}-\overline{y}_{j-1}\|_{L^2(\W_B)}\bigg(\int_{t_{j-1}}^{t_j}(t_j-s)^{-2\be}ds\bigg)^{\frac{1}{2}}.\nonumber\\
\end{align}
Note that
\begin{align}
&\int_{t_{i-1}}^{t_i}(t_{j-1}-s)^{-\a}ds\nonumber\\
=&\frac{h^{1-\a}}{1-\a}[(j-i)^{1-\a}-(j-1-i)^{1-\a}]\nonumber\\
\le&\frac{h^{1-\a}}{1-\a}[(j+1-i)^{1-\a}-(j-1-i)(j+1-i)^{-\a}]\nonumber\\
=&\frac{2h^{1-\a}}{1-\a}(j+1-i)^{-\a}\nonumber.
\end{align}
Similarly,
$$
\int_{t_{i-1}}^{t_i}(t_{j-1}-s)^{-2\be}ds\le\frac{2h^{1-2\be}}{1-2\be}(j+1-i)^{-2\be}.
$$
Therefore,
\begin{align*}
&\|\Psi_j^{h,\T}(y)-\Psi_j^{h,\T}(\overline{y})\|_{L^2(\W_B)}\\
\le&\|y_{j-1}-\overline{y}_{j-1}\|_{L^2(\W_B)}+Ch^{\g}\sum_{i=1}^{j}(j+1-i)^{-(\a\ve2\be)}\|y_{i-1}-\overline{y}_{i-1}\|_{L^2(\W_B)}.
\end{align*}
By Assumption 2.3, we have
\begin{align}\label{T4.3}
&\Big\| \sum_{j=1}^nh(t_n-(t_{j-1}+\T_jh))^{-\a}(b(\Psi_j^{h,\T}(y))-b(\Psi_j^{h,\T}(\overline{y})))\Big\|_{L^2(\W_B)}\nonumber\\
\le&\hat{L}_1\sum_{j=1}^nh(t_n-(t_{j-1}+\T_jh))^{-\a}\big\|\Psi_j^{h,\T}(y)-\Psi_j^{h,\T}(\overline{y})\big\|_{L^2(\W_B)}\nonumber\\
\le&\hat{L}_1\sum_{j=1}^nh(t_n-(t_{j-1}+\T_jh))^{-\a}\|y_{j-1}-\overline{y}_{j-1}\|_{L^2(\W_B)}\nonumber\\
+&Ch^{\g}\hat{L}_1\sum_{j=1}^nh(t_n-(t_{j-1}+\T_jh))^{-\a}\sum_{i=1}^{j}(j+1-i)^{-(\a\ve2\be)}\|y_{i-1}-\overline{y}_{i-1}\|_{L^2(\W_B)}
\end{align}
Combining this with(\ref{T5.1}), and using the independence of $\T_j$ and $y_i-\overline{y_i}, i\in\mathbb{S}^0_{j-1}$ we have
\begin{align}
&\Big\| \sum_{j=1}^nh(t_n-(t_{j-1}+\T_jh))^{-\a}(b(\Psi_j^{h,\T}(y))-b(\Psi_j^{h,\T}(\overline{y})))\Big\|_{L^2(\W)}\nonumber\\
\le&\hat{L}_1\sum_{j=1}^n\big\|h(t_n-(t_{j-1}+\T_jh))^{-\a}\big\|_{L^2(\W_\T)}\|y_{j-1}-\overline{y}_{j-1}\|_{L^2(\W)}\nonumber\\
+&Ch^{\g}\Big(n\sum_{j=1}^n\E_\T\big[|h(t_n-(t_{j-1}+\T_jh))^{-\a}|^2\big]\E_\T\Big[\sum_{i=1}^{j}(j+1-i)^{-max\{\a,2\be\}}\|y_{i-1}-\overline{y}_{i-1}\|_{L^2(\W_B)}\Big]^2\Big)^{\frac{1}{2}}\nonumber\\
\le& Ch^{1-\a}\sum_{j=1}^n(n+1-j)^{-\a }\|y_{j-1}-\overline{y}_{j-1}\|_{L^2(\W)}\nonumber\\
+&Ch^{\g+1-\a}\bigg(n\sum_{j=1}^n[(n+1-j)^{1-2\a}-(n-j)^{1-2\a}]\bigg)^{\frac{1}{2}}\Big\|\sum_{i=1}^n(n+1-i)^{-max\{\a,2\be\}}\|y_{i-1}-\overline{y}_{i-1}\|_{L^2(\W_B)}\Big\|_{L^2(\W_\T)}\nonumber\\
\le& Ch^{1-\a}\sum_{j=1}^n(n+1-j)^{-\a }\|y_{j-1}-\overline{y}_{j-1}\|_{L^2(\W)}\nonumber\\
+&Ch^{\g}T^{1-\a}\sum_{j=1}^n(n+1-j)^{-(\a\ve2\be)}\|y_{j-1}-\overline{y}_{j-1}\|_{L^2(\W)}\nonumber\\
\le& Ch^{\g}T^{1-\a}\sum_{j=1}^n(n+1-j)^{-(\a\ve2\be)}\|y_{j-1}-\overline{y}_{j-1}\|_{L^2(\W)},
\end{align}
where $C$ is a constant depending on $ \hat{L}_1, \a,\be$. Substituting this into (\ref{T3.1}) yields
\begin{align}
&\|\Phi_n^{h,\T}(y)-\Phi_n^{h,\T}(\overline{y})\|_{L^2(\W)}\nonumber\\
\le&\Big\|\sum_{j=1}^nh(t_n-(t_{j-1}+\T_jh))^{-\a}(b(\Psi_j^{h,\T}(y))-b(\Psi_j^{h,\T}(\overline{y})))\Big\|_{L^2(\W)}\nonumber\\
+&Ch^{\g}T^{\frac{1}{2}-\be}\sum_{j=1}^n(n+1-j)^{-max\{\a,2\be\}}\|y_{j-1}-\overline{y}_{j-1}\|_{L^2(\W)}\nonumber\\
\le&CTh^{\g}\sum_{j=1}^n(n+1-j)^{-max\{\a,2\be\}}\|y_{j-1}-\overline{y}_{j-1}\|_{L^2(\W)},
\end{align}
where $C$ is a positive constant depending on $ \hat{L}_1, \a,\be,\|y_j\|_{L^2(\W)}, \|\overline{y}_j\|_{L^2(\W)}, j\in\mathbb{S}^0_{n-1}$.
\eproof

In the following lemma, we use the randomized quadrature rule introduced in \cite{KW19} for integrals of stochastic processes which plays an important role in the error analysis of the randomized Milstein scheme.
\begin{lemma}
Let Assumption 2.3 hold. For each ${n\in\mathbb{S}_N}$, it holds that
\begin{equation}
\Big\|\sum_{j=1}^nh(t_n-(t_{j-1}+\T_jh))^{-\a}b(X(t_{j-1}+\T_jh))-\int_0^{t_n}(t_n-s)^{-\a}b(X(s))ds\Big\|_{L^2(\W)}
\le CTh^{\g_1}.
\end{equation}
where $\g_1=\min\{1-\a,1-\be\}$ and $C$ is a positive constant depending on $ \hat{L}_1, \a, C_0, C_1, C_2$.
\end{lemma}
\noindent
{\it Proof}.
For each $j\in\mathbb{S}_n$, due to  $\T_j\sim{\cal U}(0,1)$, we have
 \begin{equation}
 \int_{t_{j-1}}^{t_j}(t_n-s)^{-\a}b(X(s))ds=\E_\T[h(t_n-(t_{j-1}+\T_jh))^{-\a}b(X(t_{j-1}+\T_jh))].
 \end{equation}
 Further, for every $m\in\mathbb{S}_n$, set
 \begin{equation}
\left\{
\begin{array}{ll}
A_m=\sum_{j=1}^m\big(h(t_n-(t_{j-1}+\T_jh))^{-\a}b(X(t_{j-1}+\T_jh))-\int_{t_{j-1}}^{t_j}(t_n-s)^{-\a}b(X(s))ds\big)\\
A_0=0
\end{array}
\right.
\end{equation}
Then for $0\le m_1\le m_2\le n$, $m_1,m_2\in\mathbb{Z}$ it holds true that
\begin{align}
&\E_\T[(A_{m_2}-A_{m_1})|{\cal F}_{m_1}^\T]\nonumber\\
=&\sum_{j=m_1+1}^{m_2}\E_\T\Big[\Big(h(t_n-(t_{j-1}+\T_jh))^{-\a}b(X(t_{j-1}+\T_jh))-\int_{t_{j-1}}^{t_j}(t_n-s)^{-\a}b(X(s))ds\Big)\Big|{\cal F}_{m_1}^\T\Big]\nonumber\\
=&\sum_{j=m_1+1}^{m_2}\E_\T[h(t_n-(t_{j-1}+\T_jh))^{-\a}b(X(t_{j-1}+\T_jh))]-\int_{t_{m_1}}^{t_{m_2}}(t_n-s)^{-\a}b(X(s))ds=0.
\end{align}
Consequently, $(A_m)_{m\in\hat{\mathbb{S}}_n}$ is an $({\cal F}_m^\T)_{m\in\hat{\mathbb{S}}_n}$-adapted $L^p(\W_\T)$-matingale.
By an application of Theorem 2.1, we have
\begin{align}
&\|A_n\|_{L^2(\W)}^2\nonumber\\
\le&C_p^2\E_B\Big[\E_\T\Big(\sum_{j=1}^n\Big|\int_{t_{j-1}}^{t_j}[(t_n-(t_{j-1}+\T_jh))^{-\a}b(X(t_{j-1}+\T_jh))-(t_n-s)^{-\a}b(X(s))]ds\Big|^2\Big)\Big]\nonumber\\
\le&C_p^2\sum_{j=1}^n\Big\|\int_{t_{j-1}}^{t_j}[(t_n-(t_{j-1}+\T_jh))^{-\a}b(X(t_{j-1}+\T_jh))-(t_n-s)^{-\a}b(X(s))]ds\Big\|^2_{L^2(\W)}.
\end{align}
Then,
\begin{align}
&\|A_n\|_{L^2(\W)}\nonumber\\
\le&2C_p\Big(\sum_{j=1}^n\Big\|\int_{t_{j-1}}^{t_j}(t_n-s)^{-\a}(b(X(t_{j-1}+\T_jh))-b(X(s)))ds\Big\|^2_{L^2(\W)}\Big)^{\frac{1}{2}}\nonumber\\
+&2C_p\Big(\sum_{j=1}^n\Big\|\int_{t_{j-1}}^{t_j}[(t_n-(t_{j-1}+\T_jh))^{-\a}-(t_n-s)^{-\a}]b(X(t_{j-1}+\T_jh))ds\Big\|^2_{L^2(\W)}\Big)^{\frac{1}{2}}.
\end{align}
Applying Theorem 2.5 and Hölder's inequality yields
\begin{align}
&2C_p\Big(\sum_{j=1}^n\Big\|\int_{t_{j-1}}^{t_j}(t_n-s)^{-\a}(b(X(t_{j-1}+\T_jh))-b(X(s)))ds\Big\|^2_{L^2(\W)}\Big)^{\frac{1}{2}}\nonumber\\
=&2C_p\Big(\sum_{j=1}^n\E\Big|\int_{t_{j-1}}^{t_j}(t_n-s)^{-\a}(b(X(t_{j-1}+\T_jh))-b(X(s)))ds\Big|^2\Big)^{\frac{1}{2}}\nonumber\\
\le&2C_ph^{\frac{1}{2}}\Big(\sum_{j=1}^n\int_{t_{j-1}}^{t_j}(t_n-s)^{-2\a}\E|b(X(t_{j-1}+\T_jh))-b(X(s))|^2ds\Big)^{\frac{1}{2}}\nonumber\\
\le&2\hat{L}_1C_pC_2h^{\frac{1}{2}+\g}\Big(\sum_{j=1}^n\int_{t_{j-1}}^{t_j}(t_n-s)^{-2\a}ds\Big)^{\frac{1}{2}}\nonumber\\
\le&CT^{\frac{1}{2}-\a}h^{\g+\frac{1}{2}}\nonumber,
\end{align}
where $C$ is a positive constant depending on $ \hat{L}_1, \a, C_2$. Similarly,
\begin{align}
&2C_p\Big(\sum_{j=1}^n\Big\|\int_{t_{j-1}}^{t_j}[(t_n-(t_{j-1}+\T_jh))^{-\a}-(t_n-s)^{-\a}]b(X(t_{j-1}+\T_jh))ds\Big\|^2_{L^2(\W)}\Big)^{\frac{1}{2}}\nonumber\\
\le&2C_pC_0C_1h^{\frac{1}{2}}\Big(\sum_{j=1}^n\int_{t_{j-1}}^{t_j}|(t_n-s)^{-\a}-(t_n-t_{j-1})^{-\a}|^2ds\Big)^{\frac{1}{2}}\nonumber\\
\le&2C_pC_0C_1h^{\frac{1}{2}}\Big(\sum_{j=1}^{n-1}\int_{t_{j-1}}^{t_j}\Big|\a\int_{t_{j-1}}^s(t_n-r)^{-\a-1}dr\Big|^2ds+\frac{4h^{1-2\a }}{1-2\a }\Big)^{\frac{1}{2}}\nonumber\\
\le&Ch^{\frac{1}{2}}\Big(h^{1-2\a }\sum_{j=1}^{n-1}(n-j)^{-2(\a+1)}+h^{1-2\a }\Big)^{\frac{1}{2}}\nonumber\\
\le&Ch^{1-\a},
\end{align}
where $C$ is a positive constant depending on $\a, C_0, C_1$. This completes the proof.
\eproof

\begin{theorem}
Let Assumption 2.3 and Assumption 2.4 hold, then
\begin{equation}
\max_{n\in\mathbb{S}_N}\|X_n^h-X(t_n)\|_{L^2(\W)}\le Ch^{\min\{1-2\be, 1-\a\}}
\end{equation}
where $C$ is a positive constant depending on $\hat{L}_1, \hat{L}_2, \a, \be, C_0, C_1,C_2, C_3, T$.
\end{theorem}
\noindent
{\it Proof}.
For each ${n\in\mathbb{S}_N}$, it can be seen that
\begin{align}\label{T5.6}
&X_n^h-X(t_n)\nonumber\\
=&\Phi_n^{h,\T}(z)-\Phi_n^{h,\T}(\overline{z})\nonumber\\
+&\sum_{j=1}^nh(t_n-(t_{j-1}+\T_jh))^{-\a}b(X(t_{j-1}+\T_jh))-\int_0^{t_n}(t_n-s)^{-\a}b(X(s))ds\nonumber\\
+&\sum_{j=1}^nh(t_n-(t_{j-1}+\T_jh))^{-\a}(b(\Psi_j^{h,\T}(\overline{z}))-b(X(t_{j-1}+\T_jh)))\nonumber\\
+&\sum_{j=1}^n\int_{t_{j-1}}^{t_j}(t_n-s)^{-\be}\s(X(t_{j-1}))dB_s\nonumber\\
+&\sum_{j=1}^n\int_{t_{j-1}}^{t_j}(t_n-s)^{-\be}\s^{'}(X(t_{j-1}))\Big(\sum_{k=1}^{j-1}\int_{t_{k-1}}^{t_k}[(s-r)^{-\a}-(t_{j-1}-r)^{-\a}]b(X(t_{k-1}))dr\Big)dB_s\nonumber\\
+&\sum_{j=1}^n\int_{t_{j-1}}^{t_j}(t_n-s)^{-\be}\s^{'}(X(t_{j-1}))\Big(\int_{t_{j-1}}^s(s-r)^{-\a}b(X(t_{j-1}))dr\Big)dB_s\nonumber\\
+&\sum_{j=1}^n\int_{t_{j-1}}^{t_j}(t_n-s)^{-\be}\s^{'}(X(t_{j-1}))\Big(\sum_{k=1}^{j-1}\int_{t_{k-1}}^{t_k}[(s-r)^{-\be}-(t_{j-1}-r)^{-\be}]\s(X(t_{k-1}))dB_r\Big)dB_s\nonumber\\
+&\sum_{j=1}^n\int_{t_{j-1}}^{t_j}(t_n-s)^{-\be}\s^{'}(X(t_{j-1}))\Big(\int_{t_{j-1}}^s(s-r)^{-\be}\s(X(t_{j-1}))dB_r\Big)dB_s\nonumber\\
-&\int_0^{t_n}(t_n-s)^{-\be}\s(X(s))dB_s,
\end{align}
where $z=(X_0^h,\cdots,X_{n-1}^h), \overline{z}=(X(t_0),\cdots,X(t_{n-1}))$. A Taylor expansion gives
$$
\s(X(s))=\s(X(\eta(s)))+\s^{'}(X(\eta(s)))(X(s)-X(\eta(s)))+\delta(s),
$$
where $|\delta(s)|<C|X(s)-X(\eta(s))|^2$. Then
\begin{align}\label{T5.7}
&\int_0^{t_n}(t_n-s)^{-\be}\s(X(s))dB_s\nonumber\\
=&\sum_{j=1}^n\int_{t_{j-1}}^{t_j}(t_n-s)^{-\be}\s(X(t_{j-1}))dB_s\nonumber\\
+&\sum_{j=1}^n\int_{t_{j-1}}^{t_j}(t_n-s)^{-\be}\s^{'}(X(t_{j-1}))\Big(\sum_{k=1}^{j-1}\int_{t_{k-1}}^{t_k}[(s-r)^{-\a}-(t_{j-1}-r)^{-\a}]b(X(r))dr\Big)dB_s\nonumber\\
+&\sum_{j=1}^n\int_{t_{j-1}}^{t_j}(t_n-s)^{-\be}\s^{'}(X(t_{j-1}))\Big(\int_{t_{j-1}}^s(s-r)^{-\a}b(X(r))dr\Big)dB_s\nonumber\\
+&\sum_{j=1}^n\int_{t_{j-1}}^{t_j}(t_n-s)^{-\be}\s^{'}(X(t_{j-1}))\Big(\sum_{k=1}^{j-1}\int_{t_{k-1}}^{t_k}[(s-r)^{-\be}-(t_{j-1}-r)^{-\be}]\s(X(r))dB_r\Big)dB_s\nonumber\\
+&\sum_{j=1}^n\int_{t_{j-1}}^{t_j}(t_n-s)^{-\be}\s^{'}(X(t_{j-1}))\Big(\int_{t_{j-1}}^s(s-r)^{-\be}\s(X(r))dB_r\Big)dB_s\nonumber\\
+&\int_0^{t_n}(t_n-s)^{-\be}\delta(s)dB_s.
\end{align}
Substituting (\ref{T5.7}) into (\ref{T5.6}) yields
\begin{align}
&\|X_n^h-X(t_n)\|_{L^2(\W)}\nonumber\\
\le&\|\Phi_n^{h,\T}(z)-\Phi_n^{h,\T}(\overline{z})\|_{L^2(\W)}\nonumber\\
+&\Big\|\sum_{j=1}^nh(t_n-(t_{j-1}+\T_jh_))^{-\a}b(X(t_{j-1}+\T_jh))-\int_0^{t_n}(t_n-s)^{-\a}b(X(s))ds\Big\|_{L^2(\W)}\nonumber\\
+&\Big\| \sum_{j=1}^nh(t_n-(t_{j-1}+\T_jh))^{-\a}(b(\Psi_j^{h,\T}(\overline{z}))-b(X(t_{j-1}+\T_jh)))\Big\|_{L^2(\W)}\nonumber\\
+&\Big\|\sum_{j=1}^n\int_{t_{j-1}}^{t_j}(t_n-s)^{-\be}\s^{'}(X(t_{j-1}))\Big(\sum_{k=1}^{j-1}\int_{t_{k-1}}^{t_k}[(s-r)^{-\a}-(t_{j-1}-r)^{-\a}](b(X(t_{k-1}))-b(X(r)))dr\Big)dB_s\Big\|_{L^2(\W)}\nonumber\\
+&\Big\| \sum_{j=1}^n\int_{t_{j-1}}^{t_j}(t_n-s)^{-\be}\s^{'}(X(t_{j-1}))\Big(\int_{t_{j-1}}^s(s-r)^{-\a}(b(X(t_{j-1}))-b(X(r)))dr\Big)dB_s\Big\|_{L^2(\W)}\nonumber\\
+&\Big\| \sum_{j=1}^n\int_{t_{j-1}}^{t_j}(t_n-s)^{-\be}\s^{'}(X(t_{j-1}))\Big(\sum_{k=1}^{j-1}\int_{t_{k-1}}^{t_k}[(s-r)^{-\be}-(t_{j-1}-r)^{-\be}](\s(X(t_{k-1}))-\s(X(r)))dB_r\Big)dB_s\Big\|_{L^2(\W)}\nonumber\\
+&\Big\| \sum_{j=1}^n\int_{t_{j-1}}^{t_j}(t_n-s)^{-\be}\s^{'}(X(t_{j-1}))\Big(\int_{t_{j-1}}^s(s-r)^{-\be}(\s(X(t_{j-1}))-\s(X(r)))dB_r\Big)dB_s\Big\|_{L^2(\W)}\nonumber\\
+&\Big\|\int_0^{t_n}(t_n-s)^{-\be}\delta(s)dB_s\Big\|_{L^2(\W)}\nonumber\\
=:&\|\Phi_n^{h,\T}(z)-\Phi_n^{h,\T}(\overline{z})\|_{L^2(\W)}+I_1^3+I_2^3+I_3^3+I_4^3+I_5^3+I_6^3+I_7^3
\end{align}
Using Theorem 3.1 means that
\begin{equation}
\|\Phi_n^{h,\T}(z)-\Phi_n^{h,\T}(\overline{z})\|_{L^2(\W)}\le CTh^{\g}\sum_{j=1}^n(n+1-j)^{-max\{\a,2\be\}}\|X_{j-1}^h-X(t_{j-1})\|_{L^2(\W)},
\end{equation}
where $C$ is a positive constant depending on $\hat{L}_1, \a, \be,C_1,C_3$. And applying Theorem 3.2 leads to
\begin{equation}
I_1^3\le  CTh^{\g_1}.
\end{equation}
By the triangle inequality and Assumption 2.3,
\begin{align}
I_2^3&=\Big\| \sum_{j=1}^nh(t_n-(t_{j-1}+\T_jh))^{-\a}(b(\Psi_j^{h,\T}(\overline{z}))-b(X(t_{j-1}+\T_jh)))\Big\|_{L^2(\W)}\nonumber\\
&\le\hat{L}_1\Big(n\E_\T\Big[\sum_{j=1}^n|h(t_n-(t_{j-1}+\T_jh))^{-\a}|^2\E_B|\Psi_j^{h,\T}(\overline{z})-X(t_{j-1}+\T_jh))|^2\Big]\Big)^{\frac{1}{2}}
\end{align}
Note that
\begin{align}\label{T5.2}
&\|\Psi_j^{h,\T}(\overline{z})-X(t_{j-1}+\T_jh))\|_{L^2(\W_B)}\nonumber\\
\le&\bigg{\|}\int_0^{t_{j-1}}[(t_{j-1}+\T_jh-s)^{-\a}-(t_{j-1}-s)^{-\a}](b(X(\eta(s)))-b(X(s)))ds\bigg{\|}_{L^2(\W_B)}\nonumber\\
+&\bigg{\|}\int_0^{t_{j-1}}[(t_{j-1}+\T_jh-s)^{-\be}-(t_{j-1}-s)^{-\be}](\s(X(\eta(s)))-\s(X(s)))dB_s\bigg{\|}_{L^2(\W_B)}\nonumber\\
+&\bigg{\|}\int_{t_{j-1}}^{t_{j-1}+\T_jh}(t_{j-1}+\T_jh-s)^{-\a}(b(X(t_{j-1}))-b(X(s)))ds\bigg{\|}_{L^2(\W_B)}\nonumber\\
+&\bigg{\|}\int_{t_{j-1}}^{t_{j-1}+\T_jh}(t_{j-1}+\T_jh-s)^{-\be}(\s(X(t_{j-1}))-\s(X(s)))dB_s\bigg{\|}_{L^2(\W_B)}
\end{align}
Then by Assumption 2.3 and Theorem 2.5 we obtain
\begin{align}\label{T5.3}
&\bigg{\|}\int_0^{t_{j-1}}[(t_{j-1}+\T_jh-s)^{-\a}-(t_{j-1}-s)^{-\a}](b(X(\eta(s)))-b(X(s)))ds\bigg{\|}_{L^2(\W_B)}\nonumber\\
\le&\hat{L}_1\bigg(\rho_{t_j,t_{j-1}}\int_0^{t_{j-1}}|(t_{j}-s)^{-\a}-(t_{j-1}-s)^{-\a}|\E_B|X(\eta(s))-X(s)|^2ds\bigg)^{\frac{1}{2}}\nonumber\\
\le&Ch^{\g}\rho_{t_j,t_{j-1}}\le Ch^{\g+1-\a},
\end{align}
where (\ref{L4.2}) is used. Similarly, we derive
\begin{equation}
\bigg{\|}\int_{t_{j-1}}^{t_{j-1}+\T_jh}(t_{j-1}+\T_jh-s)^{-\a}(b(X(t_{j-1}))-b(X(s)))ds\bigg{\|}_{L^2(\W_B)}
\le Ch^{\g+1-\a},
\end{equation}
and
\begin{equation}
\bigg{\|}\int_{t_{j-1}}^{t_{j-1}+\T_jh}(t_{j-1}+\T_jh-s)^{-\be}(\s(X(t_{j-1}))-\s(X(s)))dB_s\bigg{\|}_{L^2(\W_B)}
\le Ch^{\g+\frac{1}{2}-\be}.
\end{equation}
Note that
\begin{align}\label{T5.4}
&\bigg{\|}\int_0^{t_{j-1}}[(t_{j-1}+\T_jh-s)^{-\be}-(t_{j-1}-s)^{-\be}](\s(X(\eta(s)))-\s(X(s)))dB_s\bigg{\|}_{L^2(\W_B)}\nonumber\\
\le&\hat{L}_1C_2h^\g\bigg(\int_0^{t_{j-1}}|(t_j-s)^{-\be}-(t_{j-1}-s)^{-\be}|^2ds\bigg)^{\frac{1}{2}}\nonumber\\
\le &Ch^{\g+\frac{1}{2}-\be}
\end{align}
Here the estimate of $\int_0^{t_{j-1}}|(t_j-s)^{-\be}-(t_{j-1}-s)^{-\be}|^2ds$ is deduced from
\begin{align}
&\int_0^{t_{j-1}}|(t_j-s)^{-\be}-(t_{j-1}-s)^{-\be}|^2ds\nonumber\\
=&(\be)^2\int_0^{t_{j-1}}\bigg(\int_{t_{j-1}}^{t_j}(r-s)^{-\be-1}dr\bigg)^2ds\nonumber\\
=&2(\be)^2\int_{t_{j-1}\le r_1<r_2\le t_j}\int_0^{t_{j-1}}(r_1-s)^{-\be-1}(r_2-s)^{-\be-1}dsdr_1dr_2\nonumber\\
\le&2(\be)^2\int_{t_{j-1}\le r_1<r_2\le t_j}\int_0^{t_{j-1}}(r_1-s)^{-\be-1}(r_2-t_{j-1})^{-\be-1}dsdr_1dr_2\nonumber\\
\le&2\be\int_{t_{j-1}\le r_1<r_2\le t_j}(r_1-t_{j-1})^{-\be}(r_2-t_{j-1})^{-\be-1}dr_1dr_2\nonumber\\
=&\frac{2\be}{(1-\be)(1-2\be)}h^{1-2\be}\nonumber.
\end{align}
Substituting (\ref{T5.3})-(\ref{T5.4}) into (\ref{T5.2}) yields
\begin{equation}\label{T5.5}
\|\Psi_j^{h,\T}(\overline{z})-X(t_{j-1}+\T_jh)\|_{L^2(\W_B)}\le Ch^{2\g}.
\end{equation}
Thus,
\begin{align}
I_2^3
&\le Ch^{2\g}\Big(n\sum_{j=1}^n\E_\T| h(t_n-(t_{j-1}+\T_jh))^{-\a}|^2\Big)^{\frac{1}{2}}\nonumber\\
&\le Ch^{2\g}\Big(\frac{nh^{2(1-\a)}}{1-2\a}\sum_{j=1}^n\big[(n+1-j)^{1-2\a}-(n-j)^{1-2\a}\big]\Big)^{\frac{1}{2}}\nonumber\\
&\le CT^{1-\a}h^{2\g}.
\end{align}
The estimate of $I_3^3$ is deduced from Assumption 2.3, Assumption 2.4 Theorem 2.5, BDG's inequality and Jensen's inequality,
\begin{align}\label{T5.6}
&I_3^3\nonumber\\
=&\Big\|\int_0^{t_n}(t_n-s)^{-\be}\s^{'}(X(\eta(s)))\Big(\int_0^{\eta(s)}[(s-r)^{-\a}-(\eta(s)-r)^{-\a}](b(X(\eta(r)))-b(X(r)))dr\Big)dB_s\Big\|_{L^2(\W)}\nonumber\\
\le&2\hat{L}_2\bigg(\int_0^{t_n}(t_n-s)^{-2\be}Q(s)ds\bigg)^{\frac{1}{2}}
\end{align}
where
$$
Q(s)=\E\Big|\int_0^{\eta(s)}[(s-r)^{-\a}-(\eta(s)-r)^{-\a}](b(X(\eta(r)))-b(X(r)))dr\Big|^2.
$$
Note that 
\clearpage
\begin{align}\label{T5.7}
&Q(s)\nonumber\\
\le&\hat{L}_1^2\int_0^{\eta(s)}|(s-r)^{-\a}-(\eta(s)-r)^{-\a}|dr
\int_0^{\eta(s)}|(s-r)^{-\a}-(\eta(s)-r)^{-\a}|\E|X(\eta(r))-X(r)|^2dr\nonumber\\
\le&\hat{L}_1^2C_1^2h^{2\g}\Big(\int_0^{\eta(s)}|(s-r)^{-\a}-(\eta(s)-r)^{-\a}|dr\Big)^2\nonumber\\
\le&Ch^{2(\g+1-\a)},
\end{align}
where the last line follows from (\ref{L4.2}). Substituting (\ref{T5.7}) into (\ref{T5.6}) yields
\begin{align}
I_3^3\le Ch^{\g+1-\a}\Big(\int_0^{t_n}(t_n-s)^{-\be}ds\Big)^{\frac{1}{2}}\le CT^{\frac{1}{2}-\be}h^{\g+1-\a}\le CT^{\frac{1}{2}-\be}h^{2\g}
\end{align}
Similarly, it is easy to verify that
\begin{align}
&I_4^3\le CT^{\frac{1}{2}-\be}h^{\g+1-\a}\le CT^{\frac{1}{2}-\be}h^{2\g},\nonumber\\
&I_5^3\le CT^{\frac{1}{2}-\be}h^{\g+\frac{1}{2}-\be}\le CT^{\frac{1}{2}-\be}h^{2\g},\nonumber\\
&I_6^3\le CT^{\frac{1}{2}-\be}h^{\g+\frac{1}{2}-\be}\le CT^{\frac{1}{2}-\be}h^{2\g}.\nonumber\\
\end{align}
Moreover, by BDG's inequality, Jensen's inequality and Theorem 2.5, we get
\begin{align}
I_7^3&\le C\bigg(\int_0^{t_n}(t_n-s)^{-2\be}\E|X(s)-X(\eta(s))|^4ds\bigg)^{\frac{1}{2}}\nonumber\\
&\le CT^{\frac{1}{2}-\be}h^{2\g}.
\end{align}
Consequently,
\begin{align*}
&\max_{n\in\mathbb{S}_N}\|X_n^h-X(t_n)\|_{L^p(\W)}\nonumber\\
\le&CTh^{\min\{1-2\be, 1-\a\}}+CTh^{\g}\sum_{j=1}^n(n+1-j)^{-(\a\ve2\be)}\max_{i\in\mathbb{S}_j}\|X_{i-1}^h-X(t_{i-1})\|_{L^p(\W)}.
\end{align*}
By Lemma 2.2
$$
\max_{n\in\mathbb{S}_N}\|X_n^h-X(t_n)\|_{L^p(\W)}\le Ch^{\min\{1-2\be, 1-\a\}},
$$
where $C$ is a positive constant depending on$\hat{L}_1, \hat{L}_2, \a, \be, C_0, C_1, C_2, C_3, T$. 

\section{Numerical experiments}
Consider the following equation:
\begin{equation} \label{7.1}
X(t) =X_0+\int_0^t(t-s)^{-\a} |\sin(X(s))|ds + \int_0^t(t-s)^{-\beta} \cos(X(s))dB_s, \quad t\in[0,T]
\end{equation}
The coefficients $b(x)=|\sin(x)|$ and $\s(x)=\cos(x)$ can be readily verified to satisfy Assumptions 2.3 and 2.4. Setting $T=1$, the numerical solution is obtained by using the randomized Milstein scheme with a minimum step size of $\hat{h}=2^{-8}$ and the expectation is approximated based on 500 independent sample paths. In simulating the singular stochastic integral $\int_{t_{j-1}}^{t_j}(t_n-s)^{-\be}\s^{'}(X_{j-1}^{h})\Big(\sum_{k=1}^{j-1}\int_{t_{k-1}}^{t_k}[(s-r)^{-\be}-(t_{j-1}-r)^{-\be}]\s(X_{k-1}^{h})dBr\Big)dB_s$, we apply the Riemann-Stieltjes integral to discretize the stochastic integral. Initially we discretize the integration interval $[t_{j-1}, t_j]$ to obtain a fixed value of $s$ for the each intergral $\int_{t_{k-1}}^{t_k}[(s-r)^{-\be}-(t_{j-1}-r)^{-\be}]\s(X_{k-1}^{h})dBr$. Subsequently, we further discretize the integration interval $[t_{k-1}, t_k]$ of this integral to obtain the value of $r$. The total integral value is then obtained by summing up all the approximations. In Figures 1 and 2, we depict the errors against the corresponding step sizes for the pairs $(\a, \be)$ are (0.3, 0.1) and (0.2, 0.3) and the slope of the line indicates the convergences which are consistent with the theoretical result of the randomized Milstein scheme.

\begin{figure}[H]
\centering
\includegraphics[height=9cm, width=15cm]{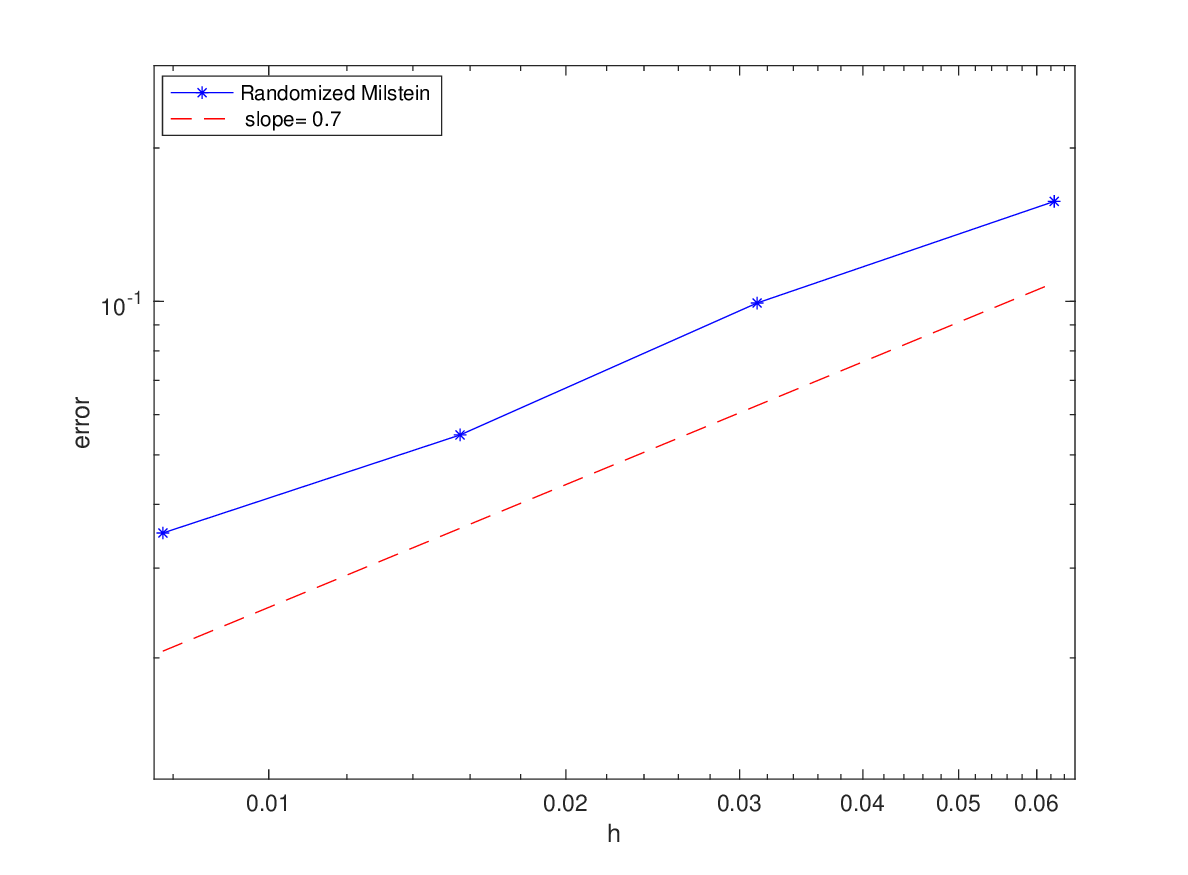}
\caption{error of the randomized Milstein scheme with $\a=0.3$, $\be=0.1$}
\label{1}
\end{figure}

\begin{figure*}[h]
\centering
\includegraphics[height=9cm, width=15cm]{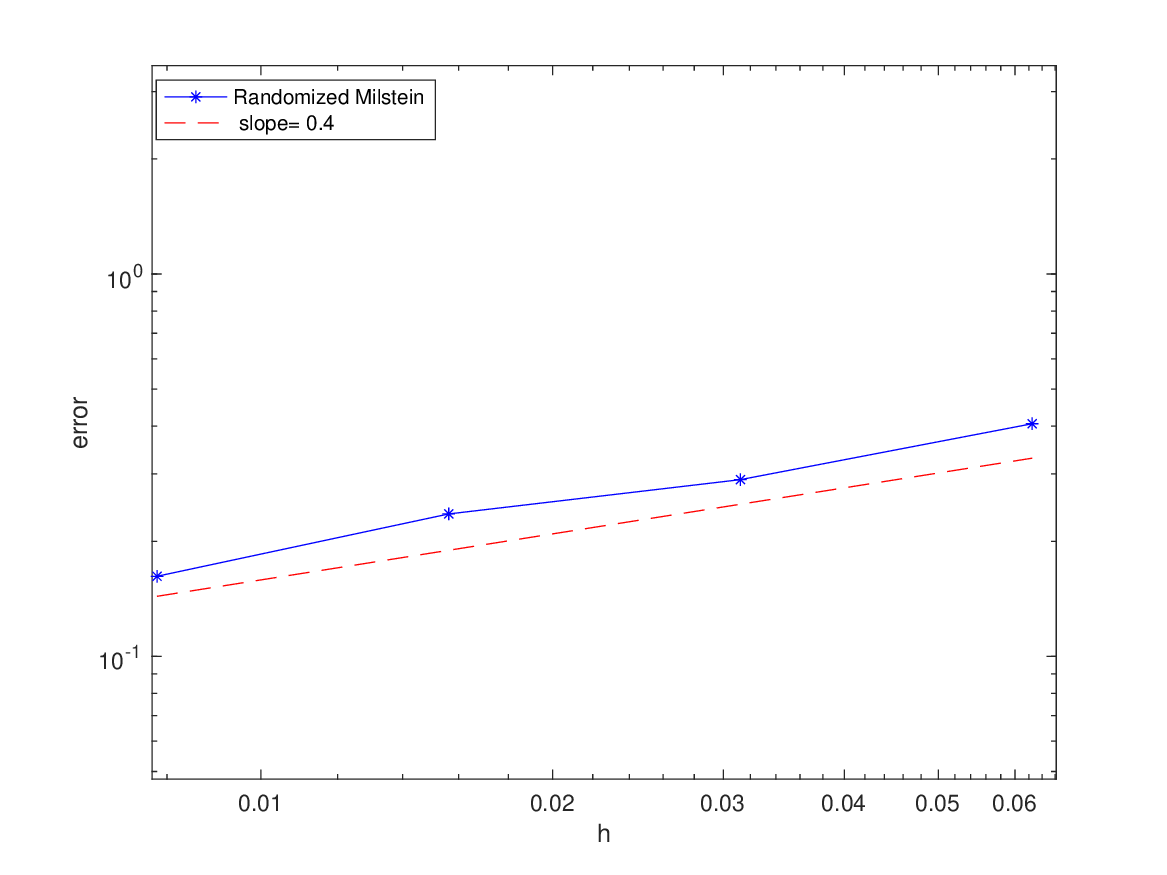}
\caption{error of the randomized Milstein scheme with $\a=0.2$, $\be=0.3$}
\label{2}
\end{figure*}

\section*{Funding}

This work is supported by the National Natural Science Foundation of China (62373383
and 62076106) and China Scholarship Council (202308310272).

\end{document}